\nonstopmode
\documentclass[leqno]{amsart}
\usepackage{amssymb}
\usepackage{graphics}

\newtheorem{theorem}{Theorem}[section]

\newtheorem{lemma}{Lemma}[section]
\theoremstyle{definition}

\newcommand{\Z}{{\mathbb Z}}

\newcommand{\sS}{{\mathcal S}}
\newcommand{\sW}{{\mathcal W}}

\newcommand{\ZZ}{{\mathbb Z}}

\begin{document}
\title{The $3x+1$ Semigroup}
\author{David Applegate}
\address{AT\&T Laboratories, Florham Park, NJ  07932-0971}
\email{david@research.att.com}

\author{Jeffrey C. Lagarias}
\address{Department of Mathematics, University of Michigan, 
Ann Arbor, MI 48109-1109}
\email{lagarias@umich.edu}
\date{April 20, 2005}

\begin{abstract}
The $3x+1$ semigroup is 
the  multiplicative semigroup $\sS$ 
of  positive rational numbers generated by
$\{ \frac{2k+1}{3k+2}: k \ge 0\}$
together with $\{2\}$. This semigroup
encodes backwards iteration under the $3x+1$ map, and the
$3x+1$ conjecture implies that it  contains every positive
integer. This semigroup is proved to be the set
of positive rationals $\frac{a}{b}$ in lowest terms
with  $b \not\equiv 0 (\bmod~3)$, 
and so contains all positive integers.
\end{abstract}

\maketitle

%
%
\section{Introduction}

The  $3x+1$
problem concerns the behavior under iteration of the $3x+1$ function
$T: \ZZ \to \ZZ$
given by
\begin{equation}\label{eq11}
T(n) = \left\{\begin{array}{ll}
\frac{n}{2} & \mbox{if } n \equiv 0 \pmod{2}\\
\frac{3n+1}{2} & \mbox{if } n \equiv 1 \pmod{2}
\end{array}\right.
\end{equation}
The $3x+1$ Conjecture asserts that for each $n \ge 1$, some iterate
$T^{(k)} (n) = 1$. This is a notoriously hard problem, work on which 
is surveyed in Lagarias \cite{La85}
and Wirsching  \cite{Wi98}.
It has been verified for all $n < 2.8 \times 10^{17}$
(see Oliveira e Silva \cite{os99} and Roosendaal \cite{Roo04}) 
but remains unsolved.

Recently H. Farkas \cite{fa04}
proposed an interesting  weakening of the $3x+1$ problem, as follows.
Let $\sS$ denote the multiplicative semigroup of positive rational
numbers generated by $\{ \frac{n}{T(n)}: n \ge 0 \}$, i.e. by 
$2$ and by $\{ \frac{2k+1}{3k+2}: k \ge 0\}$.
We call $\sS$ the $3x+1$ semigroup, and write
$$
\sS := <2, \frac{1}{2}, \frac{3}{5}, \frac{5}{8}, \frac{7}{11} \cdots >.
$$
H. Farkas  formulated the following conjecture. \\

\noindent {\bf Weak $3x+1$ Conjecture.} {\em The $3x+1$ 
semigroup $\sS$ contains every positive integer.} \\

The semigroup $\sS$ encodes inverse iteration by the $3x+1$ function.
That is, the semigroup $\sS$ contains $1 = 2 \cdot
\frac{1}{2}$, and has the property that if $T(n) \in \sS$, then also
$n \in \sS$, because each $\frac{n}{T(n)}$ is a generator of $\sS$.
 It follows that
if the $3x+1$ iteration eventually takes $n$ to $1$,
then $n$ belongs to $\sS$. Thus  the $3x+1$ conjecture
implies the weak $3x+1$ conjecture.

The  weak $3x+1$
conjecture appears  a potentially easier question to resolve
than the $3x+1$ conjecture,
since the semigroup $\sS$ permits 
some representations of integers as 
 products of generators
not corresponding to $3x+1$ iteration. Indeed, the
 object of this paper is to prove 
the following 
result characterizing  all  elements of the $3x+1$ semigroup, 
which implies the weak $3x+1$ conjecture. \\

\begin{theorem}~\label{th01}
The $3x+1$ semigroup $\sS$ equals the set of
all positive rationals $\frac{a}{b}$
in lowest terms having the property that $b \not\equiv 0 \pmod{3}$. 
In particular, it contains every positive integer. 
\end{theorem}

In order to prove this result, we shall need to study the
inverse semigroup $\sW:=\sS^{-1}$ generated by 
$\{ \frac{T(n)}{n}: n \ge 1\}$, 
i.e. by $\frac{1}{2}$
and by $\{\frac{3k+2}{2k+1}: n\ge 0\}.$ That is,
$$
\sW := \sS^{-1} = < \frac{1}{2}, \frac{2}{1}, \frac{5}{3}, 
\frac{8}{5}, \cdots>.
$$
We call this  semigroup 
the {\em wild semigroup}, following the terminology
used in a paper  \cite{La04} of the second author, which was  
inspired by the novel ``The Wild Numbers'' (\cite{Sch98}).
The paper   \cite{La04} formulated  the following conjecture. \\

\noindent {\bf Wild Numbers Conjecture.} {\em 
 The integers in the wild semigroup 
$\sW$ consist of all integers $m \ge 1$ with
$m \not\equiv 0 \pmod{3}$. Equivalently, the 
$3x+1$ semigroup $\sS$ contains
all unit fractions  $\frac{1}{m}$ such that 
$m \not\equiv 0 \pmod{3}$.} \\

Theorem~\ref{th01} is equivalent to the 
truth of both the
weak $3x+1$ conjecture and the wild numbers conjecture.
In \cite{La04} the two conjectures were shown to be
equivalent, so to deduce Theorem ~\ref{th01}
it would suffice to prove
either one of them separately.
However in  the approach taken here we 
consider them together, and  
prove them simultaneously 
using an inductive 
method in which the truth of the conjectures to given  bounds
implies their truth to a larger bound. We use a see-saw method
that increases the bound first of one, then the other.

In \S2 we show the relevance of the $3x+1$ iteration to
the weak $3x+1$ conjecture. This is the  new ingredient introduced
here relative to \cite{La04}. In \S3 we then prove properties
of integers in $\sW$ and in \S4 we complete the
argument for Theorem~\ref{th01}.

%
%

\section{Modified $3X+1$ iterations}

To prove the weak $3x+1$ conjecture by induction
on the size of the integer $m$, it
 would suffice
to prove that under forward iteration of the $3x+1$ map
starting at a given  $m \ge 2$,
 we eventually arrive at a smaller integer
$m'$, which would belong to 
the semigroup $\sS$ by the induction hypothesis.
The sequence of reverse $3x+1$ iterates going from
$m'$ back to $m$ are multiplications by elements of $\sS$,
and this would establish that  $m \in \sS.$ 
However, if this argument could
be carried out, it would prove more,
namely the $3x+1$ conjecture
itself. Since this problem  seems out of reach, we 
considered a modification of  this approach.

We take advantage of the fact that the $3x+1$ iteration
decreases ``almost all'' integers, in the sense of
\cite[Theorem A]{La85}. We recall that
forward iteration of the $3x+1$ function $T(\cdot)$ for $j$ steps
is known to decrease the value of an integer $n$ in
most congruence classes $n \pmod{2^j}$. Recall that the
first $j$ steps of the $3x+1$ iteration
are uniquely determined by the class
$n \pmod{2^j}$ and that every symbol pattern of even and odd integers
of length $j$ occurs in some trajectory
of length $j$, cf. Lagarias\cite[Theorem B]{La85}.
A residue class $ s \pmod{2^j}$ is said
to have  a {\em strong stopping time} $k \le j$ if 
the smallest integer $s \ge 2$ in the residue
class decreases after $k$ steps of iteration.
This property is then inherited by all members $\ge 2$ of
the residue class.
As $j$ increases the fraction of integers not having
a strong stopping time 
goes to zero, but there still remain exponentially many
residue classes $\pmod{2^j}$ not having  this decreasing property
(\cite[Theorems C and D ]{La85}).

The semigroup $\sS$ permits
the possibility of going ``uphill'' by taking an initial
value $n$ to a value $mn$ via some integer multiplier $m$,
provided  $m \in \sS^{-1}= \sW$. That is, if  $\frac{1}{m} \in \sS$
and if we know   $mn \in \sS$ then we may deduce
 $n = \frac{1}{m} \cdot mn \in \sS.$
We pay a price in going ``uphill''
 of increasing the initial size of the integer,
but in doing so we may move from a ``bad'' residue 
class $ s \pmod{2^j}$ to a ``good'' residue class 
$ ms \pmod{2^j}$ which under iteration results in 
such a large decrease in the size of the number  
that it overcomes  the added multiplicative factor $m$ and arrives
at an integer smaller than $n$ in $\le j$ steps. 
One can use this procedure  only for $j$ steps ahead because 
the members of the residue class only possess the same symbolic
dynamics for $j$ steps, and we wish the property of decrease
to hold for all members of the residue class.
If so, we  
 can carry out the induction step for all integers in
this particular  ``bad'' residue 
class $ s \pmod{2^j}$. Another variation of this idea
is to multiply by various $m$'s in the middle of the
first $j$ steps of the iteration; there is no reason why
the multiplication must be done only at the first step, one
may still gain by switching the residue class in the
middle of the iteration.

One can now ask: is there a finite $j$ and 
a finite list $\{ m_1, m_2, \cdots m_r\}$ 
of integer elements in $\sW$
such that suitable multiplications by elements of this
list will decrease elements in  every residue class
$\pmod{2^j}$ in this fashion? If so, this would yield
a proof of the weak $3x+1$ conjecture by induction on $n$.

 This approach comes very 
close to succeeding, but 
there is an obstruction that in principle
prevents it from succeeding.
We found by computer search, for small values of $j$,
multiplier lists  that 
established decrease for  every residue
class $\pmod{2^j}$ except for the class $-1 \pmod{2^j}$.
These searches revealed that the class  $-1 \pmod{2^j}$
resisted elimination for $12 \le j \le 30$.
We then looked for and found the following proof
that the class  $-1 \pmod{2^j}$  can never 
be eliminated by this method.
The iterates of a positive integer $n$ 
in the congruence class  $-1 \pmod{2^j}$   will  behave 
exactly the
same way as $-1$ does for the first $j$ steps, 
allowing multipliers. We may write  the
$j$-th iterate of $-1$ obtained using multipliers as
$\frac{m_1m_2\cdots m_j a (-1) + b}{2^j}$,
in which   $m_k$ is the multiplier used at the $k$-th  step
(we allow $m_k=1$),
$a$ is a power of $3$,  and $b$ is a positive integer.
For this multiplier sequence 
any $n \equiv -1 \pmod{2^j}$ will map to 
$\frac{m_1m_2 \cdots m_j a n + b}{2^j}$ after $j$
steps. However we 
must have 
$$
\frac{(m_1m_2\cdots m_j) a (-1) + b}{2^j} \le -1,
$$
because all iterates of $-1$, times multipliers, remain
negative. Rearranging this inequality gives  
$$
(m_1m_2\cdots m_j) a \ge 2^j + b.
$$
Now, for positive $n$,
 multiplying both sides by $\frac{n}{2^j}$ yields 
$$
\frac{(m_1m_2\cdots m_j) a n + b}{2^j} \ge n + \frac{b(n+1)}{2^{j}} > n.
$$ 
It follows  that  decrease
cannot have  occurred after $j$ steps, and
an  argument for no decrease
at any intermediate step is similar. 

We conclude that  to get an inductive proof of
the weak $3x+1$ conjecture along these lines,
a new method will be needed to handle integers in the 
``bad'' congruence class $-1 \pmod{2^j}$, and it will be necessary
to consider an infinite set of multipliers in $\sW$.

We now prove the decrease mentioned above for 
all residue classes $ \pmod{4096}$ except  
the class $-1 \pmod{4096}$, 
using a fixed finite set $H$ of multipliers given below;
these are residue classes $\pmod{2^j}$ for $j=12$.
In what follows  it will be important
that the decrease is by a constant factor strictly 
smaller than one. In \S3 we will verify the
hypothesis  $H \subset \sW$ made  in this lemma. 
%
%
%
%

\begin{lemma}
\label{lem:4096}
If $H = \{5, 7, 11, 13, 23, 29, 43\} \subset \sW$, then
for every integer $x > 1$ with $x \not\equiv -1 \pmod{4096}$ 
there exists $s \in
\sW$ such that $sx \in \Z$ and $sx \le \frac{76}{79} x$.
\end{lemma}

%
%
%

\begin{table}[htbp]
\center
\renewcommand{\arraystretch}{0.7}
\begin{tabular}{|l|r|r|l|}
\hline
class & \multicolumn{1}{|c|}{asymptotic} & \multicolumn{1}{|c|}{worst-case} & path\\
class bits & \multicolumn{1}{|c|}{ratio} & \multicolumn{1}{|c|}{ratio} & \\
\hline
$0 \pmod{2}$ & $1/2$ & $1/2^*$ & $2 \rightarrow 1$\\
0 & 0.5000 & 0.5000 & \\
\hline
$1 \pmod{4}$ & $3/4$ & $4/5^*$ & $5 \rightarrow 8 \rightarrow 4$\\
10 & 0.7500 & 0.8000 & \\
\hline
$3 \pmod{16}$ & $9/16$ & $2/3$ & $3 \rightarrow 5 \rightarrow 8 \rightarrow 4 \rightarrow 2$\\
1100 & 0.5625 & 0.6667 & \\
\hline
$11 \pmod{32}$ & $27/32$ & $10/11$ & $11 \rightarrow 17 \rightarrow 26 \rightarrow 13 \rightarrow 20 \rightarrow 10$\\
11010 & 0.8438 & 0.9091 & \\
\hline
$27 \pmod{128}$ & $117/128$ & $25/27$ & $27 \rightarrow 41 * 13 = 533 \rightarrow 800 \rightarrow 400 \rightarrow $\\
1101100 & 0.9141 & 0.9259 & $\;\;\;\rightarrow 200 \rightarrow 100 \rightarrow 50 \rightarrow 25 $\\
\hline
$91 \pmod{256}$ & $225/256$ & $80/91$ & $91 * 25 = 2275 \rightarrow 3413 \rightarrow 5120 \rightarrow $\\
11011010 & 0.8789 & 0.8791 & $\;\;\;\rightarrow 2560 \rightarrow 1280 \rightarrow 640 \rightarrow 320 \rightarrow $\\
&&& $\;\;\;\rightarrow 160 \rightarrow 80 $\\
\hline
$219 \pmod{256}$ & $243/256$ & $209/219$ & $219 \rightarrow 329 \rightarrow 494 \rightarrow 247 \rightarrow 371 \rightarrow $\\
11011011 & 0.9492 & 0.9543 & $\;\;\;\rightarrow 557 \rightarrow 836 \rightarrow 418 \rightarrow 209 $\\
\hline
$59 \pmod{128}$ & $81/128$ & $38/59$ & $59 \rightarrow 89 \rightarrow 134 \rightarrow 67 \rightarrow 101 \rightarrow $\\
1101110 & 0.6328 & 0.6441 & $\;\;\;\rightarrow 152 \rightarrow 76 \rightarrow 38 $\\
\hline
$123 \pmod{256}$ & $189/256$ & $91/123$ & $123 * 7 = 861 \rightarrow 1292 \rightarrow 646 \rightarrow $\\ 
11011110 & 0.7383 & 0.7398 & $\;\;\;\rightarrow 323 \rightarrow 485 \rightarrow 728 \rightarrow 364 \rightarrow $\\
&&& $\;\;\;\rightarrow 182 \rightarrow 91 $\\
\hline
$251 \pmod{256}$ & $207/256$ & $203/251$ & $251 * 23 = 5773 \rightarrow 8660 \rightarrow 4330 \rightarrow $\\
11011111 & 0.8086 & 0.8088 & $\;\;\;\rightarrow 2165 \rightarrow 3248 \rightarrow 1624 \rightarrow 812 \rightarrow $\\
&&& $\;\;\;\rightarrow 406 \rightarrow 203 $\\
\hline
$7 \pmod{64}$ & $45/64$ & $5/7$ & $7 * 5 = 35 \rightarrow 53 \rightarrow 80 \rightarrow 40 \rightarrow 20 \rightarrow $\\
111000 & 0.7031 & 0.7143 & $\;\;\;\rightarrow 10 \rightarrow 5$\\
\hline
$39 \pmod{128}$ & $105/128$ & $32/39$ & $39 * 35 = 1365 \rightarrow 2048 \rightarrow 1024 \rightarrow $\\
1110010 & 0.8203 & 0.8205 & $\;\;\;\rightarrow 512 \rightarrow 256 \rightarrow 128 \rightarrow 64 \rightarrow 32 $\\
\hline
$103 \pmod{512}$ & $351/512$ & $71/103$ & $103 \rightarrow 155 \rightarrow 233 * 13 = 3029 \rightarrow $\\
111001100 & 0.6855 & 0.6893 & $\;\;\;\rightarrow 4544 \rightarrow 2272 \rightarrow 1136 \rightarrow 568 \rightarrow $\\
&&& $\;\;\;\rightarrow 284 \rightarrow 142 \rightarrow 71$\\
\hline
$359 \pmod{512}$ & $315/512$ & $221/359$ & $359 * 35 = 12565 \rightarrow 18848 \rightarrow 9424 \rightarrow $\\
111001101 & 0.6152 & 0.6156 & $\;\;\;\rightarrow 4712 \rightarrow 2356 \rightarrow 1178 \rightarrow 589 \rightarrow $\\
&&& $\;\;\;\rightarrow 884 \rightarrow 442 \rightarrow 221$\\
\hline
$231 \pmod{256}$ & $135/256$ & $122/231$ & $231 * 5 = 1155 \rightarrow 1733 \rightarrow 2600 \rightarrow $\\
11100111 & 0.5273 & 0.5281 & $\;\;\;\rightarrow 1300 \rightarrow 650 \rightarrow 325 \rightarrow 488 \rightarrow $\\
&&& $\;\;\;\rightarrow 244 \rightarrow 122 $\\
\hline
$23 \pmod{32}$ & $27/32$ & $20/23$ & $23 \rightarrow 35 \rightarrow 53 \rightarrow 80 \rightarrow 40 \rightarrow 20$\\
11101 & 0.8438 & 0.8696 & \\
\hline
\end{tabular}
\caption{Decreasing weak $3x+1$ paths, for $x \not\equiv 15 \pmod{16}$).}
\label{tab1a}
\end{table}

%
%
%
%

\begin{table}[htbp]
\center
\renewcommand{\arraystretch}{0.7}
\begin{tabular}{|l|r|r|l|}
\hline
class & \multicolumn{1}{|c|}{asymptotic} & \multicolumn{1}{|c|}{worst-case} & path\\
class bits & \multicolumn{1}{|c|}{ratio} & \multicolumn{1}{|c|}{ratio} & \\
\hline
$15 \pmod{128}$ & $81/128$ & $10/15$ & $15 \rightarrow 23 \rightarrow 35 \rightarrow 53 \rightarrow 80 \rightarrow 40 \rightarrow $\\
1111000 & 0.6328 & 0.6667 & $\;\;\;\rightarrow 20 \rightarrow 10 $\\
\hline
$79 \pmod{256}$ & $243/256$ & $76/79$ & $79 \rightarrow 119 \rightarrow 179 \rightarrow 269 \rightarrow 404 \rightarrow $\\
11110010 & 0.9492 & 0.9620 & $\;\;\;\rightarrow 202 \rightarrow 101 \rightarrow 152 \rightarrow 76 $\\
\hline
$207 \pmod{256}$ & $225/256$ & $182/207$ & \small $207 * 5 = 1035 \rightarrow 1553 * 5 = 7765 \rightarrow $\\
11110011 & 0.8789 & 0.8792 & \small $\;\;\;\rightarrow 11648 \rightarrow 5824 \rightarrow 2912 \rightarrow $\\
&&& \small $\;\;\;\rightarrow 1456 \rightarrow 728 \rightarrow 364 \rightarrow 182 $\\
\hline
$47 \pmod{128}$ & $117/128$ & $43/47$ & $47 * 13 = 611 \rightarrow 917 \rightarrow 1376 \rightarrow $\\
1111010 & 0.9141 & 0.9149 & $\;\;\;\rightarrow 688 \rightarrow 344 \rightarrow 172 \rightarrow 86 \rightarrow 43 $\\
\hline
$111 \pmod{128}$ & $99/128$ & $86/111$ & \small $111 * 11 = 1221 \rightarrow 1832 \rightarrow 916 \rightarrow $\\
1111011 & 0.7734 & 0.7748 & \small $\;\;\;\rightarrow 458 \rightarrow 229 \rightarrow 344 \rightarrow 172 \rightarrow 86 $\\
\hline
$31 \pmod{64}$ & $33/64$ & $16/31$ & $31 * 11 = 341 \rightarrow 512 \rightarrow 256 \rightarrow $\\
111110 & 0.5156 & 0.5161 & $\;\;\;\rightarrow 128 \rightarrow 64 \rightarrow 32 \rightarrow 16 $\\
\hline
$63 \pmod{128}$ & $99/128$ & $49/63$ & $63 * 11 = 693 \rightarrow 1040 \rightarrow 520 \rightarrow $\\
1111110 & 0.7734 & 0.7778 & $\;\;\;\rightarrow 260 \rightarrow 130 \rightarrow 65 \rightarrow 98 \rightarrow 49 $\\
\hline
$127 \pmod{256}$ & $129/256$ & $64/127$ & $127 * 43 = 5461 \rightarrow 8192 \rightarrow 4096 \rightarrow $\\
11111110 & 0.5039 & 0.5039 & $\;\;\;\rightarrow 2048 \rightarrow 1024 \rightarrow 512 \rightarrow 256 \rightarrow $\\
&&& $\;\;\;\rightarrow 128 \rightarrow 64 $\\
\hline
$255 \pmod{512}$ & $387/512$ & $193/255$ & \small $255 * 43 = 10965 \rightarrow 16648 \rightarrow 8224 \rightarrow $\\
111111110 & 0.7559 & 0.7569 & \small $\;\;\;\rightarrow 4112 \rightarrow 2056 \rightarrow 1028 \rightarrow $\\
&&& \small $\;\;\;\rightarrow 514 \rightarrow 257 \rightarrow 386 \rightarrow 193$\\
\hline
$511 \pmod{1024}$ & $783/1024$ & $391/511$ & \small $511 \rightarrow 767 * 29 = 22243 \rightarrow 33365 \rightarrow $\\
1111111110 & 0.7646 & 0.7652 & \small $\;\;\;\rightarrow 50048 \rightarrow 25024 \rightarrow 12512 \rightarrow $\\
&&& \small $\;\;\;\rightarrow 6256 \rightarrow 3128 \rightarrow 1564 \rightarrow $\\
&&& \small $\;\;\;\rightarrow 782 \rightarrow 391 $\\
\hline
$1023 \pmod{2048}$ & $1089/2048$ & $544/1023$ & \small $1023 * 11 = 11253 \rightarrow 16880 \rightarrow $\\
11111111110 & 0.5317 & 0.5318 & \small $\;\;\;\rightarrow 8440 \rightarrow 4220 \rightarrow 2110 \rightarrow $\\
&&& \small $\;\;\;\rightarrow 1055 * 11 = 11605 \rightarrow 17408 \rightarrow $\\
&&& \small $\;\;\;\rightarrow 8704 \rightarrow 4352 \rightarrow 2176 \rightarrow $\\
&&& \small $\;\;\;\rightarrow 1088 \rightarrow 544 $\\
\hline
$2047 \pmod{4096}$ & $3267/4096$ & $1633/2047$ & \small $2047 * 11 = 22517 \rightarrow 33776 \rightarrow $\\
111111111110 & 0.7976 & 0.7978 & \small $\;\;\;\rightarrow 16888 \rightarrow 8444 \rightarrow 4222 \rightarrow $\\
&&& \small $\;\;\;\rightarrow 2111 * 11 = 23221 \rightarrow 34832 \rightarrow $\\
&&& \small $\;\;\;\rightarrow 17416 \rightarrow 8708 \rightarrow 4354 \rightarrow $\\
&&& \small $\;\;\;\rightarrow 2177 \rightarrow 3266 \rightarrow 1633 $\\
\hline
\end{tabular}
\caption{Decreasing weak $3x+1$ paths, for $x \equiv 15 \pmod{16}$.}
\label{tab1b}
\end{table}

\begin{proof}
This is established
case by case in
 Tables \ref{tab1a} and \ref{tab1b}.  
Every path shown consists of iterations of
$T(\cdot)$ and multiplications by integers in $H$, and thus consists
of iterations of multiplications by elements of $\sW$. 
The iteration takes $k$ steps, where $k$ is the given number of bits,
and for integers $n$ in  the class $s \pmod{2^k}$, one has 
$T^{(k)}(n) = c(s)n + d(s)$, with 
$$
c(s) = \frac{3^{l} m_1m_2 \cdots m_k}{2^k},
$$ 
in which  the $m_i$ are the multipliers at each step and $l$ is the number
of odd elements in the resulting trajectory, and $d(s) \ge 0$. The quantity
$c(s)$ is the ``asymptotic ratio'' reported in the second column of
the tables. 

The ``class bits'' presented
 in these tables are binary strings
comprising the binary expansion
of the residue class written in reverse order. The set
of these  binary strings
together form
 a prefix code which by inspection certifies that every 
residue class $(\bmod~ 4096)$
is covered except $-1~ (\bmod~ 4096).$
The data on the far right in the table 
gives the action on the smallest positive
element in the congruence class (resp. second smallest element for
the class containing $n=1$). In each case the factor of decrease on
all elements of the progression (excluding the
element  $n=1$),
reported as the ``worst-case ratio'' in the table,
 is that given by the decrease on
this particular element. 
\end{proof}

To  deal with the residue class $-1 \pmod{2^j}$,
we next show that there always  exists a 
simple (but infinite) sequence of multipliers having
the property that, starting from $n \equiv -1 \pmod{2^j}$, with $n >0$
one arrives at a final integer $n'$ that is only slightly larger
than the initial starting point $n$. We will later make use
of this to eliminate the congruence class $-1 + 2^j \pmod{2^{j+1}}$,
in an induction on $j$.

%
%
%

\begin{lemma}
\label{lem:tjmx}
Let $x$, $k$, and $j$ be positive integers such that $x \equiv -1
\pmod{2^k}$, with $1 \le j \le k$ and 
$j \equiv 1,5 \pmod{6}$.  Then the multiplier 
$m = \frac{2^j+1}{3}$ is an integer satisfying $m \equiv 1,5 \pmod{6}$,
with  the property that the $j$-th iterate of $m x$ satisfies the bound
$$
T^j(m x) = x +
\frac{x+1}{2^j} \le \frac{2^j + 2}{2^j} x \;,
$$
and $T^j(m x) \equiv -1 \pmod{2^{k-j}}$.
If in addition $x \not\equiv -1 \pmod{2^{k+1}}$, 
then $T^j(m x) \not\equiv -1 \pmod{2^{k+1-j}}$.
\end{lemma}


\begin{proof}
Since $j \equiv 1,5 \pmod{6}$, $2^j \equiv 2,5 \pmod{9}$, 
so $m = \frac{2^j+1}{3}$
is an odd integer and $m \not\equiv 0 \pmod{3}$.  Since $mx \equiv -m
\pmod{2^k}$, $mx$ is odd, so
$$T(mx) = \frac{3mx + 1}{2} = \frac{2^j x + x + 1}{2}\;\;.$$
Since $x \equiv -1 \pmod{2^k}$ and $k \ge j$, $2^j x + x + 1 \equiv 0
\pmod{2^j}$.  Thus
$$T^j(mx) = \frac{2^j x + x + 1}{2^j} = x + \frac{x+1}{2^j}\;\;.$$
$\frac{x+1}{2^j} \equiv 0 \pmod{2^{k-j}}$, so $T^j(mx) \equiv -1
\pmod{2^{k-j}}$.  If in addition $x \not\equiv -1 \pmod{2^{k+1}}$,
then $\frac{x+1}{2^j} \not\equiv 0 \pmod{2^{k+1-j}}$, so $T^j(mx)
\not\equiv -1 \pmod{2^{k+1-j}}$.
\end{proof}

To make use of  Lemma~\ref{lem:tjmx}
in an inductive proof, we need
to establish  that after using it on $n$ to obtain
$n'=T^{(j)}(n)$ a
single  $3x+1$ iteration applied to 
$n'$ produces an integer  $n''$
smaller than $n$. This is the aim of the
following lemma, which gives an inductive  method  of eliminating the
class $-1 + 2^j \pmod{2^{j+1}}$ using a suitable integer multiplier $m$ ,
assuming that $m$ is a wild integer.

%
%
%
%
\begin{lemma}
\label{lem:onestep}
Suppose $H \subset \sW$.
Let $x \equiv -1 \pmod{2^k}$ and $x \not\equiv -1 \pmod{2^{k+1}}$,
for a fixed $k \ge 12$. Now choose $j$ 
so that $j \equiv 1 \pmod{6}$ and $k-10 \le j \le k-5$. Then 
$m = (2^j + 1)/3$ is an integer, and if  $m \in
\sW$, then there exists $s \in \sW$ such
that $sx \in \Z$ and $sx \le \frac{1235}{1264} x$. 
\end{lemma}

\begin{proof}
First, note that 
for all $k \ge 12$, since $j \equiv 1 \pmod{6}$ and $j \geq k-10 \geq 2$,
we have $j \ge 7$.
From Lemma \ref{lem:tjmx}, $m \in \Z$ and 
$T^j(mx) = x + \frac{x+1}{2^j}$, so
there exists $s_1 \in \sW$ such that $s_1 m x \in \Z$, 
$s_1 m x = x+ \frac{x+1}{2^j}$, $s_1 m x \equiv -1 \pmod{2^{k-j}}$, 
and $s_1 m x\not\equiv -1 \pmod{2^{k+1-j}}$.  
But $k-j \le 10$, so
$s_1 m x \not\equiv -1 \pmod{2^{11}}$.
Thus from Lemma \ref{lem:4096}, there exists 
$s_2 \in \sW$ such that $s_2 (s_1 m x) \in \Z$ and 
$s_2 (s_1 m x) \le 76/79(s_1 m x)$. Now 
$j \ge 7$ and the bound of 
Lemma~\ref{lem:tjmx} gives
$$x + \frac{x+1}{2^j} \le
\frac{2^j+2}{2^j} x \le \frac{130}{128} x,
$$ 
so that 
$s_2 s_1 m x \le \frac{76}{79} \frac{130}{128} x = \frac{1235}{1264} x$.
\end{proof}

%
%

\section{Wild Integers}

The wild integers are the 
integers in the wild semigroup $\sW$.
The ``multiplier'' approach begun in \S2 required the use of  
multipliers that are wild integers, and  indicated
that in taking  this approach one would need to consider 
an infinite set of multipliers. 
This in turn seems to require  understanding the complete
structure of the integer elements in $\sW$, which leads
to investigation of the wild numbers conjecture.

In this section we establish properties of
wild integers, giving  criteria for establishing
their existence.  
We first show that the elements in $H$ in \S2
are wild integers. Here we write
$g(n) = \frac{3n+2}{2n+1}$.

%
%
%
%

\begin{lemma}\label{lem31}
The set $H =  \{5, 7, 11, 13, 23, 29, 43\}$
is contained in the wild semigroup $\sW = \sS^{-1}$.
\end{lemma}
\begin{proof}
Table~\ref{table31} below gives certificates
showing that
the elements in $H$ belong to $\sW$,
representing them 
in terms of the generators of $\sW$.
The table uses the notation  $g(n) = \frac{3n+2}{2n+1}$, for $n \ge 1$.
Aside frome $p=5$,
these  identities were  found by computer search by Allan Wilks, 
see \S2 of  \cite{La04}.
\end{proof}
%
%
%
%

\begin{table}[htbp]
\begin{tabular}{|rrl|}
\hline
  $ 5$ & $ =$ & $(\frac{1}{2})^2 \cdot 
(\frac{11}{7})^2 \cdot \frac{17}{11} \cdot \frac{26}{17} 
\cdot \frac{83}{55} \cdot \frac{98}{65} \cdot \frac{125}{83} $  \\
& $ =$ & $(\frac{1}{2})^2 \cdot 
g(3)^2 \cdot g(5)\cdot g(8) \cdot g(27) \cdot g(32) \cdot g(41) $  \\

  $ 7$ & $ =$ &  $(\frac{1}{2})^2 \cdot 
 \frac{11}{7}    \cdot \frac{26}{17}    \cdot \frac{35}{23} 
\cdot \frac{215}{143} 
\cdot \frac{299}{199} \cdot  \frac{323}{215} 
\cdot \frac{371}{247} \cdot  \frac{398}{265}    $ \\
& $ =$ &  $(\frac{1}{2})^2 \cdot 
g(3) \cdot g(8) \cdot g(11) \cdot g(71) \cdot g(99) \cdot g(107) 
\cdot g(123) \cdot g(132)$ \\

  $ 11$ & $ = $ & $(\frac{1}{2})^2 \cdot
(\frac{11}{7}) ^2 \cdot \frac{26}{17} \cdot \frac{35}{23}  \cdot 
\frac{215}{143} \cdot \frac{299}{199} \cdot \frac{323}{215} 
\cdot  \frac{371}{247}    \cdot \frac{398}{265} $ \\
& $ = $ & $(\frac{1}{2})^2 \cdot
g(3)^2 \cdot g(8) \cdot g(11)  \cdot g(71) \cdot g(99) \cdot g(107)
\cdot g(123) \cdot g(132)$ \\

  $ 13$  & $ =$ & $(\frac{1}{2})^3 \cdot
(\frac{11}{7})^2 \cdot (\frac{17}{11})^3 \cdot (\frac{26}{17})^2 
\cdot \frac{35}{23} \cdot 
\frac{215}{143} \cdot \frac{299}{199} \cdot 
\frac{323}{215}    \cdot  \frac{371}{247}     \cdot \frac{398}{265}  $   \\
& $ =$ & $(\frac{1}{2})^3 \cdot
g(3)^2 \cdot g(5)^3 \cdot g(8)^2 \cdot g(11) \cdot g(71) \cdot g(99) \cdot
g(107)\cdot g(123) \cdot g(132)$   \\

  $  23$  & $ =$ & $(\frac{1}{2})^5 \cdot
 \frac{11}{7}\cdot \frac{26}{17} \cdot \frac{35}{23} 
\cdot \frac{47}{31} \cdot \frac{137}{91} \cdot  \frac{155}{103} \cdot 
\frac{206}{137} \cdot \frac{215}{143} \cdot (\frac{299}{199})^2
\cdot \frac{323}{215} \cdot \frac{353}{235} \cdot 
\frac{371}{247} ~\cdot $\\
 & & ~~~~~$\cdot~ (\frac{398}{265})^2 \cdot \frac{530}{353}$ \\
& $ =$ & $(\frac{1}{2})^5 \cdot
g(3) \cdot g(8) \cdot g(11)  
\cdot g(15) \cdot g(45) \cdot  g(51) \cdot g(68) \cdot g(71) \cdot g(99)^2
~\cdot $ \\
 & & ~~~~~$\cdot~ g(107) \cdot g(117) \cdot g(123) \cdot g(132)^2 \cdot g(176)$ \\

  $  29$  & $ =$ & $(\frac{1}{2})^5 \cdot
(\frac{11}{7})^4 \cdot (\frac{17}{11})^2 \cdot (\frac{26}{17})^2 
\cdot \frac{29}{19} 
\cdot \frac{38}{25} 
\cdot  (\frac{83}{55})^2 \cdot (\frac{98}{65})^2 \cdot (\frac{125}{83})^2 $ \\
& $ =$ & $(\frac{1}{2})^5 \cdot
g(3)^4 \cdot g(5)^2 \cdot g(8)^2 \cdot g(9) \cdot g(12) 
\cdot g(27)^2 \cdot g(32)^2 \cdot g(41)^2 $ \\

 $  43$ & $ =$ & $(\frac{1}{2})^{11} \cdot
 (\frac{11}{7})^5 \cdot (\frac{17}{11})^2 \cdot(\frac{26}{17})^3 
\cdot \frac{29}{19} \cdot \frac{35}{23} \cdot \frac{38}{25} \cdot
(\frac{83}{55})^2 \cdot (\frac{98}{65})^2 \cdot (\frac{125}{87})^2 
\cdot \frac{215}{143} ~\cdot $  \\
& &~~~~~  $ \cdot~ \frac{299}{199}    \cdot \frac{305}{203} \cdot 
\frac{323}{215} \cdot \frac{344}{229} \cdot \frac{371}{247} 
\cdot \frac{398}{265} \cdot \frac{458}{305}$ \\
& $ =$ & $(\frac{1}{2})^{11} \cdot
 g(3)^5 \cdot g(5)^2 \cdot g(8)^3 \cdot g(9) \cdot g(11) \cdot g(12) \cdot
g(27)^2 \cdot g(32)^2 \cdot g(41)^2 ~\cdot $  \\
& &~~~~~  $ \cdot~  g(71) \cdot g(99) \cdot g(101) \cdot g(107) \cdot g(114) \cdot g(123) 
\cdot g(132) \cdot g(152)$ \\
\hline
\end{tabular}
\caption{Membership certificates in $\sW$ for members of $H$.}
\label{table31}
\end{table}

The following lemma uses the truth of the 
 weak $3x+1$ conjecture on an initial interval
to extend the range on which the wild numbers
conjecture holds.

%
%

\begin{lemma}\label{lem32}
Suppose that the weak $3x+1$ conjecture holds for 
$1 \le n \le 2^j -2$
and that the wild numbers conjecture holds for
$1 \le m \le \frac{2^{j}-1}{189}$, with $j \ge 16$.
Then the wild numbers  conjecture holds for 
$1 \le m \le \frac{2^{j+1}-1}{189}.$
\end{lemma}

\begin{proof}
It suffices to prove that every prime $q$ with
$\frac{2^j-1}{189} < q \leq \frac{2^{j+1}-1}{189}$ lies in
$\sW$. Proceeding by induction on increasing $q$, 
we may assume every prime $p$
with $3 < p < q$ lies in
$\sW$. It now suffices to prove: 
{\em There exists a positive integer $n \le 2^j-2$ with
 $nq \in \sW$.}
For if so, then the induction hypothesis  implies 
that $n \in \sS$ so
that $q = \frac{1}{n}\cdot nq \in \sW$. In establishing
this we will consider only those
$n$ such that $nq \equiv -1$ (mod 9). 
Then $nq = 3l + 2$ for
some positive integer $l$, and $nq = t\cdot (2l+1)$ where $t =
\frac{3l+2}{2l+1} \in \sW$. Thus it will  suffice to show
$2l+1 \in \sW$. \\

To carry this out, define
 $a$ as the least positive residue with  $aq \equiv -1$ (mod 9),
so that $0 < a < 9$. For $n$ in the arithmetic progression
$n = 9k + a$, setting $nq=3l+2$, we have 
$$
2l+1 = 2(\frac{nq-2}{3})+1 =
\frac{2}{3}((9k+a)q-2)+1 = 6qk + r, ~~\mbox{with}~~ r := \frac{2aq-1}{3}.
$$
The condition  $aq \equiv -1$ (mod 9) gives 
 $r \equiv -1 ~(\bmod~6)$, and  $r~(\bmod~6q)$ is invertible $(\bmod~6q)$.
For the given prime $q$ the values $a$ and $r$ are
determined,  and we need to find a suitable
value of $k$. If $0 \le k < 6q$ then: 
$$
n = 9k + a\leq 9(6q-1) + a < 54q \leq 54(\frac{2^{j+1}-1}{189}) =
\frac{2}{7}(2^{j+1}-1) \le 2^j-2,
$$
so $n \in \sS$ by hypothesis. 
Therefore it suffices to prove: 
{\em For each prime $q$ with
$\frac{2^j-1}{189} < q \leq \frac{2^{j+1}-1}{189}$ there
 exists an integer $0 \le k < 6q$ such that
$6qk + r \in \sW$.}

Define a positive integer to be $q$-smooth if all its prime
factors are smaller than $q$. Let $\Sigma_q$ denote the set of
$q$-smooth integers $s$ with $0 < s < 6q$ and $gcd(s, 6q) = 1$.
Then every $s \in \Sigma_q$ is a product of primes $p$ with 
$5 \le p < q$, and the induction
hypothesis implies that $s \in \sW$.\\

{\em Claim. If $q \geq 256$ then $|\Sigma_q| > q - 1$.}

Assuming the 
claim is true, we can apply it in our
situation because $q > \frac{2^j-1}{189} \geq \frac{2^{16}-1}{189}
> 346$. The claim implies that $\Sigma_q$
contains more than half of the invertible residue classes (mod
$6q$), since $\phi(q) = 2(q-1)$.
 Therefore in the group of invertible residue classes (mod
$6q$), the sets $\Sigma_q$ and $r\cdot\Sigma_q^{-1}$ must meet,
since each contains more than half of the classes. Therefore $s_1
\equiv r\cdot s_2^{-1}$ (mod $6q$), for some $s_1,
s_2 \in \Sigma_q$. Now $s_1 s_2 \equiv r$ (mod $6q$), and we may 
define $k \ge 0$ by setting $s_1s_2 = 6qk + r$. Since each $s_i \in
\Sigma_q \subseteq \sW$ we have $6qk + r \in
\sW$. Since $s_1, s_2 < 6q$ we find that $k < 6q$, as required.
Thus the  proof of Lemma 3.2 will be complete once the
claim is established.

To prove the claim, since $\phi(6q)=2q-2$
we may reformulate it as the assertion: {\it  there are
at most $q-2$ invertible residue classes below $6q$ which
are not $q$-smooth.} The non-$q$-smooth numbers below $6q$ relatively
prime to $q$ consist of the primes $q'$ with $q < q' < 6q$
together with  the integers  $5q'$ where $q'$ is prime
with $q < q' < \frac{6}{5} q.$ Thus we
must show that for $q > 256$, 
\begin{equation}~\label{307}
(\pi(6q) - \pi(q)) + \pi(\frac{6}{5}q) - \pi(q) \le q-2.
\end{equation}
The left side of (\ref{307}) is $O(\frac{q}{\log q})$
by the prime number theorem, so (\ref{307})  holds for all
sufficiently large $q$; it remains to
establish the specific bound.  We use explicit inequalities for prime
counting functions due to 
Rosser and Schoenfeld \cite[Theorems 1 and 2]{RS62},
which state that for all $x \ge 17$, 
\begin{equation}\label{N303}
\frac{x}{\log x} < \pi(x) < \frac{x}{\log x - \frac{3}{2}},
\end{equation}
and also that, for all $x \ge 114$,
$$
\pi(x) < \frac{5}{4}  \frac{x}{\log x}.
$$
The first of these inequalities gives  
$$
\pi(6x) \le \frac{6x}{ \log(6x) - \frac{3}{2}} \le \frac{6x}{\log x}
$$
since $\log 6 \ge \frac{3}{2}.$ The second gives, 
for $x \ge 256$, 
\begin{eqnarray*}
\pi (\frac{6}{5}x) 
& < & \frac{5}{4}  
\left(\frac{\frac{6}{5}x}{\log (\frac{6}{5}x)}\right) \\
 & < & \frac{3}{2}  \frac{x}{\log x} 
\left(\frac{\log x}{\log x + \log \frac{6}{5}}\right) \\
 & < &  \frac{3}{2}  \frac{x}{\log x} \left( 1 - 
\frac{\frac{1}{6}}{\log x + \frac{1}{6}} \right) \\
& < & \frac{3}{2}  \frac{x}{\log x} - 2, 
\end{eqnarray*}
where we used $\log \frac{6}{5} > \frac{1}{6}$, and $x \ge 256$
was used at the last step. Combining these bounds gives,
for $x \ge 256 > e^{11/2}$, 
$$
\pi(6x) + \pi (\frac{6}{5}x) - 2 \pi(x)  < 
\frac{11}{2} \frac{x}{\log x} - 2 \le x - 2,
$$
which proves the claim.
\end{proof}
%
%
\section{Completion of Proofs}

\begin{proof}
[Proof of Theorem \ref{th01}]
The theorem is equivalent to the  truth of the weak $3x+1$ conjecture
and the wild numbers conjecture. Together these two 
conjectures imply
that the semigroup $\sS$ contains all rationals 
$\frac{a}{b} = a \cdot \frac{1}{b}$ with $b \not\equiv 0
\pmod{3}$. However $\sS$  contains no rational 
$\frac{a}{b}$ in lowest terms with $b \equiv 0 \pmod{3}$,
because no generator of $\sS$ contains a multiple of $3$ in
its denominator. Conversely, if $\sS$ contains all such
rationals, then both conjectures hold.

We prove the weak $3x+1$ conjecture and wild numbers conjecture
simultaneously by induction on $k \ge 12$, using the following
three inductive hypotheses.

(1) For each  integer $x >1$ with  $x \not\equiv -1  \pmod{2^k}$ 
there is an element $s \in \sW$ such that $sx$ is an integer
and 
$sx \le  \frac{1235}{1264}x$. 

(2)  The weak $3x+1$ conjecture is true for $1 \le n \le 2^{k}-2.$

(3) The wild integers conjecture is true for 
$1 \le m \le \frac{2^k-1}{189}.$

We treat the induction step first, and the base case afterwards.
We suppose the inductive hypotheses hold for some $k \ge 12$,
and must show show they then hold for $k+1$.

Hypotheses  (2) and (3) for $k$ permit 
Lemma~\ref{lem32} to apply,  
whence  for $k \ge 16$ we conclude that   
inductive hypothesis (3) holds  
for $k+1$. For the remaining
cases $12 \le k \le 15$ we verify 
inductive hypothesis (3) for $k+1$ 
directly by computation, which is included
 in the base case  below.

Inductive hypothesis (1) 
for $k$ gives  that all elements smaller than  $2^{k+1}-1$
except possibly $2^{k}-1$ can be decreased by multiplication by
an element of $\sW$ to a smaller integer.
We wish to  apply Lemma~\ref{lem:onestep} to show that
all elements in the congruence class
$-1 + 2^k  \pmod{2^{k+1}}$ can also be decreased by multiplication
by an element of $\sW$ to an integer smaller 
by the multiplicative factor $\frac{1235}{1264}$.
First,  Lemma~\ref{lem31} shows that the elements of $H$ belong to $\sW$,
establishing one hypothesis of Lemma~\ref{lem:onestep}.
Second,  the other  multiplier 
$m = \frac{2^{j}+1}{3}$ in the hypothesis of Lemma \ref{lem:onestep}
has  
$j = k - 5 - \left(k \pmod{6}\right)$, and satisfies
$m \le \frac{2^{(k+1)-6} +1}{3} \le \frac{2^{k+1}-1}{189}$, 
so $m \in \sW$ 
by inductive hypothesis (3), which is already 
established to hold for $k+1$. Thus 
all the hypotheses of Lemma~\ref{lem:onestep} are satisfied,
and its conclusion verifies 
the inductive hypothesis (1) for $k+1$.

Next,  inductive hypothesis (1) for $k+1$  establishes
the decreasing property for all integers 
$1 < n \le 2^{k+1} -2$, hence the weak $3x+1$ conjecture
follows for all integers in this range. This verifies
inductive hypothesis (2) for $k+1$, and so  
completes the induction step.

It remains to treat the base case, which is $k=12$
for hypotheses (1) and (2), 
and $k=16$ for hypothesis (3). 
For $k=12$ inductive hypothesis (1) is verified
by Lemma~\ref{lem:4096}. The inductive  hypothesis (2)
for $k=12$ is verified by the fact that the
$3x+1$ conjecture has been  checked over the range 
$1 \le n \le 2^{12} = 4096.$

Finally we must  verify  inductive  hypothesis (3) 
for $k=16$. This requires 
verifying the wild numbers conjecture for
$1 \le m \le \frac{2^{16}-1}{189} = \frac{65535}{189} < 400$.
It suffices to do this for all primes below $400$, except $p=3$.
Representations  in the generators of
$\sW$ for all such  primes below $50$, are given in \cite{La04}.
(Table 3 gives representations for some of these primes.)
For primes $50 \le p \le 400$, one can check
the criterion by computer using the method of Lemma~\ref{lem32},
finding  by
computer search  a $q$-smooth number in the appropriate arithmetic
progression, for each prime $q$ in the interval,
and using the truth of the $3x+1$ conjecture for 
$1 < x < 10^{5}$. In fact, this $q$-smooth calculation can 
be carried out by computer for every $q$ with $ 11 < q < 400$,
and only the certificates  for $p=5$, $p=7$, $p=11$ in Table~\ref{table31}
are needed to begin the induction. As an example,
for $q=13$ we have $a=2$ and $r= \frac{2aq - 1}{3}= 17$,
and the arithmetic progression $78k + 17$ contains $875=5^3 \cdot 7$.
\end{proof}

%
%
\section{Concluding Remarks}

The proofs in this paper are computer-intensive.
Computer experimentation 
played an important role in the discovery of the
patterns underlying the induction. This included the
efficacy of using multipliers to eliminate congruence
classes $\pmod{2^k}$ in \S2, and in uncovering
the existence of
the ``intractable'' residue class $-1 \pmod{2^k}$.
If one had studied the problem without using the computer,
the ``intractable''  case $-1 \pmod{2^k}$ could have been 
uncovered first, and this  might have discouraged further
investigation of this proof approach.
It was also important to have the evidence detailed
in \cite{La04}, which provided a strong element of
confidence in the truth of the weak $3x+1$ conjecture
and wild numbers conjecture. 

Extensive computations were needed  to
find the data in the tables. Once  found,
this data in the tables provides  
``succinct certificates'' for checking correctness
of the congruence class properties, 
which can be 
verified by hand. Similarly the induction step
is in principle checkable by hand.

The proof methods developed in
this paper should apply more generally in determining the  integers in
various  multiplicative semigroups of rationals
having a similar nature. 

\paragraph{\bf Acknowledgments.} The authors thank Allan Wilks
for helpful computations concerning the wild semigroup
in Lemma ~\ref{lem31}. They  thank the reviewer for
many useful comments, including a suggested revision of the proof
of Lemma~\ref{lem32} incorporated in this version.  
The second author did most of his work on this paper while
employed  at AT\&T Labs. 
%
%

\ifx\undefined\bysame
\newcommand{\bysame}{\leavevmode\hbox to3em{\hrulefill}\,}
\fi


\begin{thebibliography}{10}

\bibitem{fa04}
H. Farkas,
{\em Variants of the $3N+1$ problem and multiplicative semigroups,}
in: {\it Geometry, Spectral Theory, Groups and Dynamics:
Proceedings in Memory of Robert Brooks}
(M. Entov, Y. Pinchover and M. Sageev, Eds.),
Israel Mathematical Conference Proceedings, 
Amer. Math. Soc., Providence 2005, to appear.



\bibitem{La85}
J.~C. Lagarias, {\em The $3x+1$ problem and its generalizations}, Amer. Math.
  Monthly {\bf 92} (1985), 3--21.

\bibitem{La04}
J.~C. Lagarias,
{\em Wild and Wooley numbers}, Amer. Math. Monthly, to appear.

\bibitem{os99}
T. Oliveira e Silva,
{\em Maximum excursion and stopping time record-holders for the $3x+1$
problem: computational results,} 
Math. Comp. {\bf 68}, No. 1 (1999), 371--384.

\bibitem{RS62}
J. B. Rosser and L. Schoenfeld,
{\em Approximate formulas for some functions of prime numbers,}
Illinois J. Math. {\bf 6} (1962), 64--94.

\bibitem{Roo04}
E. Roosendaal, 
On the $3x+1$ problem, web document:
{\tt http://personal.computrain.nl/eric/wondrous/}.
(Distributed computation for $3x+1$ problem statistics.)


\bibitem{Sch98}
P. Schogt,
{\em The Wild Numbers}, Four Walls Eight Windows:
New York 1998.


\bibitem{Wi98}
G.~J. Wirsching, 
{\em The dynamical system on the natural numbers generated by
  the $3n+1$ function}, 
Lecture Notes in Math. No. 1681, Springer-Verlag:
Berlin 1998.

\end{thebibliography}
\end{document}